\documentclass{article}
\usepackage{amsmath}
\usepackage{amsfonts}

\input{tcilatex}

\begin{document}

\title{Accurate distribution of $\mathbf{X}^{\text{T}}\mathbf{X}$ with
singular, idempotent variance-covariance matrix $\mathbb{V}$}
\date{}
\author{Hao Yuan Zhang and Jan \ Vrbik \\
Department of Mathematics\\
Brock University, 500 GLenridge Ave.\\
St. Catharines, Ontario, Canada, L2S 3A1}
\maketitle

\begin{abstract}
Assume that $\mathbf{X}$ is a set of sample statistics which follow a
special case Central Limit Theorem, namely: as the sample size $n$ increases
the corresponding distribution becomes multivariate Normal with the mean (of
each $X$) equal to zero and with an \emph{idempotent} variance-covariance
matrix $\mathbb{V}$. It is well known that $\mathbf{X}^{\text{T}}\mathbf{X}$
has (in the same limit), a $\chi ^{2}$ distribution with degrees of freedom
equal to the trace of $\mathbb{V}.$ In this article we extend the above
result to include the corresponding $\frac{1}{n}$-proportional corrections,
making the new approximation substantially more accurate and extending its
range of applicability to small-size samples.
\end{abstract}

\section{Introduction}

Consider a random independent sample of size $n$ from a distribution,
resulting in $p$ sample statistics, say $X_{1},$ $X_{2},...X_{p}$
(collectively denoted $\mathbf{X}$). We assume that these are defined in
such a way (visualize each being a function of standardized sample means)
that their joint distribution converges (as $n\rightarrow \infty $) to
multivariate Normal. This implies that their joint $j^{th}$-order cumulants
can be expanded in inverse powers of $n$ as follows:%
\begin{equation*}
\mathcal{K}^{(j)}=\left( \kappa ^{(j,0)}+\frac{\kappa ^{(j,1)}}{n}+\frac{%
\kappa ^{(j,2)}}{n^{2}}+...\right) n^{1-j/2}
\end{equation*}%
where $\mathcal{K}^{(j)}$ and each of the corresponding $\kappa ^{(j,\ell )}$
is a fully symmetric \emph{tensor} with $j$ (implicit) indices.

When $j=1,$ $\mathcal{K}^{(1)}$ is a (column) \emph{vector} of the expected
values of $\mathbf{X}$; the $\kappa ^{(1,0)}$ term must always (rather
exceptionally) equal to zero, and the $\kappa ^{(j,1)}$ term we rename $%
\mathbf{\mu }$, so that%
\begin{equation*}
\mathcal{K}^{(1)}=\frac{\mathbf{\mu }}{\sqrt{n}}+...
\end{equation*}%
Similarly, $\mathcal{K}^{(2)}$ is the corresponding variance-covariance 
\emph{matrix}, expanded (and notationally simplified) as follows:%
\begin{equation}
\mathcal{K}^{(2)}=\mathbb{V}^{(0)}+\frac{\mathbb{V}^{(1)}}{n}+...  \label{K2}
\end{equation}

We will now assume that $\mathbb{V}^{(0)}$ is \emph{idempotent} with a \emph{%
trace} equal to $k.$ Our task is to find the approximate distribution of 
\begin{equation}
T=\mathbf{X}^{\text{T}}\mathbf{X}  \label{main}
\end{equation}%
to the $\frac{1}{n}$ accuracy, extending the familiar result which states
that, in the $n\rightarrow \infty $ limit, the distribution of $T$ becomes $%
\chi ^{2}$ with $k$ degrees of freedom.

\section{MGF and PDF of rotated $\mathbf{X}$}

To find a more accurate approximation for $T$, we recall that there is an
ortho-normal matrix $\mathbb{R}$ which diagonalizes $\mathbb{V}^{(0)}$ thus:%
\begin{equation*}
\mathbb{R~V}^{(0)}\mathbb{R}^{\text{T}}\equiv \mathbb{H}
\end{equation*}%
where $\mathbb{H}$ is \emph{main-diagonal}, with the first $k$ diagonal
elements equal to $1,$ the rest of them equal to $0;$ (\ref{main}) can then
be rewritten as%
\begin{equation*}
T=\mathbf{Z}^{\text{T}}\mathbf{Z}
\end{equation*}%
where%
\begin{equation*}
\mathbf{Z}\equiv \mathbb{R}~\mathbf{X}
\end{equation*}%
Note that transforming a $j^{\text{th}}$-order cumulant of $\mathbf{X}$ into
its $\mathbf{Z}$ counterpart is achieved by%
\begin{equation*}
\mathcal{\tilde{K}}_{i_{1},i_{2}...i_{j}}=\sum_{\ell _{1},\ell _{2}...\ell
_{j}=1}^{p}\mathbb{R}_{i_{1},\ell _{1}}\mathbb{R}_{i_{2},\ell _{2}}...%
\mathbb{R}_{i_{j},\ell _{j}}\mathcal{K}_{\ell _{1},\ell _{2}...\ell _{j}}
\end{equation*}%
Also note that the $\kappa ^{(j,m)}$ (individually) transform in the same
manner, and (essential to the rest of this article) that any component of
the new $\tilde{\kappa}^{(j,0)}$ with a lower (\emph{implicit}, in this
notation) index in the $k+1$ to $p$ range must be equal to zero (no longer
true for $\tilde{\kappa}^{(j,1)}$, $\tilde{\kappa}^{(j,2)}...$, including $%
\mathbf{\tilde{\mu}}$ and $\mathbb{\tilde{V}}^{(1)}$).

This implies that, in the $n\rightarrow \infty $ limit, the $Z_{i}$s with $%
i\leq k$ are independent, standardized Normal, and the remaining $Z_{i}$s ($%
p-k$ of them) are identically equal to zero (implying the $T\epsilon \chi
_{k}^{2}$ result). But again, this happens only in the $n\rightarrow \infty $
limit; both $\frac{\tilde{\mu}_{i}}{\sqrt{n}}$ and the corresponding
components of $\frac{\mathbb{\tilde{V}}^{(1)}}{n}$ may remain \emph{non-zero}
even when $i>k$ (similar to what happens to $\frac{g(\bar{X})-g(\mu )}{%
\sigma |g^{\prime }(\mu )|}$ in the same limit, where $\bar{X}$ is a sample
mean and $g$ is an arbitrary function). The $\chi _{k}^{2}$ approximation
thus becomes less accurate with decreasing $n$. To improve its accuracy, we
proceed to find a $\frac{1}{n}$-proportional correction to it.

\subsection{Expanding MGF}

The corresponding cumulant generating function of $\mathbf{Z}$, expanded to
the $\frac{1}{n}$ accuracy, is given by%
\begin{eqnarray*}
&&K(t_{1},t_{2},...t_{p})\overset{}{\equiv }K(\mathbf{t})\overset{}{=} \\
&&\frac{\mathbf{t}^{T}\mathbb{H~}\mathbf{t}}{2}+\frac{\mathbf{t}^{\text{T}}%
\mathbf{\tilde{\mu}}^{(1)}}{\sqrt{n}}+\frac{\tilde{\kappa}^{(3,0)}\circ 
\mathbf{t}\circ \mathbf{t}\circ \mathbf{t}}{6\sqrt{n}}+\frac{\mathbf{t}^{T}%
\mathbb{\tilde{V}}^{(1)}\mathbf{t}}{2n}+\frac{\tilde{\kappa}^{(4,0)}\circ 
\mathbf{t}\circ \mathbf{t}\circ \mathbf{t}\circ \mathbf{t}}{24n}+...
\end{eqnarray*}%
where the last ellipsis indicates terms beyond the $\frac{1}{n}$ accuracy,
and $\tilde{\kappa}\circ \mathbf{t}$ implies contracting the last (implicit)
index of $\tilde{\kappa}$ with the only index of $\mathbf{t}$, i.e.,
explicitly,%
\begin{equation*}
\sum_{i_{\ell }=1}^{p}\tilde{\kappa}_{i_{1},i_{2},...i_{\ell }}\cdot
t_{i_{\ell }}
\end{equation*}
Thus, for example, $\tilde{\kappa}^{(3,0)}\circ \mathbf{t}\circ \mathbf{t}%
\circ \mathbf{t}$ means that all $3$ indices of $\tilde{\kappa}^{(3,0)}$
have been contracted, one by one, with an index of $\mathbf{t}$, resulting
in a scalar.

It is relatively easy to convert the above cumulant generating function into
the corresponding moment generating function, thus 
\begin{eqnarray}
&&M(\mathbf{t})\overset{}{=}\exp [K(\mathbf{t})]\overset{}{=}\exp \left( 
\frac{\sum_{i=1}^{k}t_{i}^{2}}{2}\right) \cdot  \label{jv1} \\
&&\left( 1+\frac{\mathbf{t}^{T}\mathbb{\tilde{V}}^{(1)}\mathbf{t}}{2n}+\frac{%
\tilde{\kappa}^{(4,0)}\circ \mathbf{t}\circ \mathbf{t}\circ \mathbf{t}\circ 
\mathbf{t}}{24n}+\frac{\left( \mathbf{t}^{\text{T}}\mathbf{\tilde{\mu}}^{(1)}%
\overset{}{+}\frac{1}{6}\tilde{\kappa}^{(3,0)}\circ \mathbf{t}\circ \mathbf{t%
}\circ \mathbf{t}\right) ^{2}}{2n}+...\right)  \notag
\end{eqnarray}%
where this time we have discarded not only the $o(\frac{1}{n})$ terms, but
also terms with odd powers of $\mathbf{t}$; they do not contribute to our
final answer, as we will show shortly.

\subsection{Corresponding PDF}

The moment generating function (\ref{jv1}) converts to the following joint
PDF of the $\mathbf{Z}$ distribution:%
\begin{eqnarray}
&&f(z_{1},z_{2},...z_{p})\overset{}{=}  \label{pdf} \\
&&(2\pi )^{-k/2}\cdot \left( 1+\frac{\mathbf{D}^{\text{T}}\mathbb{\tilde{V}}%
^{(1)}\mathbf{D}}{2n}+\frac{\tilde{\kappa}^{(4,0)}\circ \mathbf{D}\circ 
\mathbf{D}\circ \mathbf{D}\circ \mathbf{D}}{24n}+\right.  \notag \\
&&\left. \frac{\left( \mathbf{D}^{\text{T}}\mathbf{\tilde{\mu}}^{(1)}+\frac{1%
}{6}\tilde{\kappa}^{(3,0)}\circ \mathbf{D}\circ \mathbf{D}\circ \mathbf{D}%
\right) ^{2}}{2n}+...\right) \exp \left( -\frac{\sum_{i=1}^{k}z_{i}^{2}}{2}%
\right) \prod\limits_{i=k+1}^{p}\delta (z_{i})  \notag
\end{eqnarray}%
where $\delta (z_{i})$ stands for the Dirac delta function (visualize it as
a PDF of a Normal distribution with $\mu =0$ and $\sigma \rightarrow 0$),
and $\mathbf{D}$ is a differential operator whose individual components are
partial derivatives with respect to $z_{i}$ (made more explicit later on).

\section{Finding CDF of $T$}

Let us now compute%
\begin{equation}
\Pr (T<u)=\Pr \left( \sum_{i=1}^{p}Z_{i}^{2}<u\right)   \label{Pr}
\end{equation}%
This can be done by integrating (\ref{pdf}) over a sphere of radius $\sqrt{u}
$ in the $p$-dimensional space of the $z_{i}$ variables. The main
contribution is of course from 
\begin{eqnarray}
&&(2\pi )^{-k/2}\idotsint\limits_{\mathcal{R}_{p}(\sqrt{u})}\exp \left( -%
\frac{\sum_{i=1}^{k}z_{i}^{2}}{2}\right) \prod\limits_{i=k+1}^{p}\delta
(z_{i})~dz_{1}dz_{2}...dz_{p}  \notag \\
&=&(2\pi )^{-k/2}\idotsint\limits_{\mathcal{R}_{k}(\sqrt{u})}\exp \left( -%
\frac{\sum_{i=1}^{k}z_{i}^{2}}{2}\right) dz_{1}dz_{2}...dz_{k}  \label{aux}
\end{eqnarray}%
where $\mathcal{R}_{d}(\rho )$ denotes a $d$-dimensional sphere of radius $%
\rho .$ Introducing $r\equiv \sqrt{\sum_{i=1}^{k}z_{i}^{2}}$ and utilizing
Fisher's geometrical method of multidimensional integration, (\ref{aux})
equals to%
\begin{eqnarray}
&&F(u)\equiv (2\pi )^{-k/2}\int_{0}^{\sqrt{u}}S_{k}(r)\exp \left( -\frac{%
r^{2}}{2}\right) dr  \notag \\
&=&(2\pi )^{-k/2}\frac{2\pi ^{k/2}}{\Gamma \left( \frac{k}{2}\right) }%
\int_{0}^{\sqrt{u}}r^{k-1}\exp \left( -\frac{r^{2}}{2}\right) dr  \notag \\
&=&\frac{1}{2^{k/2}\Gamma \left( \frac{k}{2}\right) }\int_{0}^{u}w^{k/2-1}%
\exp \left( -\frac{w}{2}\right) dw  \label{F}
\end{eqnarray}%
where%
\begin{equation*}
S_{k}(r)=\frac{2\pi ^{k/2}r^{k-1}}{\Gamma \left( \frac{k}{2}\right) }
\end{equation*}%
is the area of the surface of a $k$-dimensional sphere of radius $r,$ and (%
\ref{F}) is clearly the CDF of the $\chi _{k}^{2}$ distribution.

\subsection{$\mathbf{\tilde{\protect\mu}}^{(1)}$ and $\mathbb{\tilde{V}}%
^{(1)}$ corrections}

From what follows it is easy to see why terms of (\ref{pdf}) having an odd
power of $D_{i}$ (for any one or more $i$) must yield zero contribution to (%
\ref{Pr}) - that goes for those we have already deleted, as well for most of
those which are still explicitly a part of (\ref{pdf}), such as, for example%
\begin{equation*}
...+\frac{D_{1}\mathbb{\tilde{V}}_{12}^{(1)}D_{2}}{n}+...
\end{equation*}%
etc. The only terms which remain are those containing $D_{i}^{2},$ $%
D_{i}^{4},$ $D_{i}^{6},$ $D_{i}^{2}D_{j}^{2},$ $D_{i}^{4}D_{j}^{2},$ and $%
D_{i}^{2}D_{j}^{2}D_{\ell }^{2},$ where $i\neq j\neq k.$ Let us establish
their contribution to (\ref{Pr}). Starting with $D_{i}^{2},$ we get (using $%
D_{1}^{2}$ as a proxy; due to the obvious symmetry, all $D_{i}^{2},$ where $%
1\leq i\leq k,$ will contribute equally):%
\begin{eqnarray}
&&(2\pi )^{-k/2}\didotsint_{\mathcal{R}_{p}(\sqrt{u})}\frac{\partial ^{2}}{%
\partial z_{1}^{2}}\exp \left( -\frac{\sum_{i=1}^{k}z_{i}^{2}}{2}\right)
\dprod_{i=k+1}^{p}\delta (z_{i})~dz_{1}dz_{2}...dz_{p}  \notag \\
&=&(2\pi )^{-k/2}\didotsint_{\mathcal{R}_{k}(\sqrt{u})}\frac{\partial ^{2}}{%
\partial z_{1}^{2}}\exp \left( -\frac{\sum_{i=1}^{k}z_{i}^{2}}{2}\right)
dz_{1}dz_{2}...dz_{k}  \notag \\
&=&(2\pi )^{-k/2}\dint\limits_{-\sqrt{u}}^{\sqrt{u}}\frac{\partial ^{2}\exp
(-\frac{z_{1}^{2}}{2})}{\partial z_{1}^{2}}\didotsint_{\mathcal{R}_{k-1}(%
\sqrt{u-z_{1}^{2}})}\exp \left( -\frac{\sum_{i=2}^{k}z_{i}^{2}}{2}\right)
dz_{2}...dz_{k}~dz_{1}  \notag \\
&=&(2\pi )^{-k/2}\dint\limits_{-\sqrt{u}}^{\sqrt{u}}\frac{\partial ^{2}\exp
(-\frac{z_{1}^{2}}{2})}{\partial z_{1}^{2}}\int_{0}^{\sqrt{u-z_{1}^{2}}%
}S_{k-1}(r)\cdot \exp (-\tfrac{r^{2}}{2})dr~dz_{1}  \notag \\
&=&-\frac{2\pi ^{-1/2}}{2^{k/2}\Gamma (\frac{k-1}{2})}\dint\limits_{-\sqrt{u}%
}^{\sqrt{u}}\frac{\partial \exp (-\frac{z_{1}^{2}}{2})}{\partial z_{1}}\cdot 
\frac{\partial }{\partial z_{1}}\int_{0}^{\sqrt{u-z_{1}^{2}}}r^{k-2}\exp (-%
\tfrac{r^{2}}{2})dr~dz_{1}  \notag \\
&=&-\frac{2\pi ^{-1/2}}{2^{k/2}\Gamma (\frac{k-1}{2})}\exp (-\frac{u}{2}%
)\dint\limits_{-\sqrt{u}}^{\sqrt{u}}z_{1}^{2}(u-z_{1}^{2})^{(k-3)/2}dz_{1} 
\notag \\
&=&-\frac{2}{k\cdot 2^{k/2}\Gamma (\frac{k}{2})}\cdot u^{k/2}\exp (-\frac{u}{%
2})=-\frac{2u}{k}\cdot \frac{u^{k/2-1}\exp (-\frac{u}{2})}{2^{k/2}\Gamma (%
\frac{k}{2})}  \label{H2}
\end{eqnarray}%
To follow the derivation, one must be able to: (i) handle Dirac's delta
function inside an integral, (ii) establish the intersection of a plane with
a sphere in $k$ dimensions, (iii) perform by-part integration, (iv)
differentiate with respect to a parameter which appears in the upper limit
of an integral, and (v) be familiar with integrals relating to Beta function.

Careful examination of (\ref{pdf}) reveals that (\ref{H2}) needs to be added
to (\ref{F}) after being multiplied by 
\begin{equation}
a\equiv \frac{\sum_{i=1}^{k}\left( \mathbb{\tilde{V}}_{i,i}^{(1)}+\tilde{\mu}%
_{i}^{(1)}\cdot \tilde{\mu}_{i}^{(1)}\right) }{2n}  \label{a}
\end{equation}

As for the $D_{i}^{2}$ contribution to (\ref{Pr}) when $i>k$ (note that for
these values of $i$ there is no contribution from the $\tilde{\kappa}%
^{(3,0)} $ and $\tilde{\kappa}^{(4,0)}$ terms): it is fairly obvious that
adding the sum of squares of these $p-k$ random variables with zero variance
and the mean of $\frac{\mathbf{\mu }_{i}^{(1)}}{\sqrt{n}}$ will simply
increase the value of $T$ by%
\begin{equation*}
\frac{\sum_{i=k+1}^{p}\tilde{\mu}_{i}^{(1)}\cdot \tilde{\mu}_{i}^{(1)}}{n}
\end{equation*}%
Since the contribution of $\mathbb{\tilde{V}}_{i,i}^{(1)}$ is the same as
that of $\mu _{i}^{(1)}\cdot \mu _{i}^{(1)}$ at the MGF level, it must be
the same (to this level of approximation) at the PDF level; this implies
that the only effect of the $Z_{i}$ ($i>k$) variables will be adding%
\begin{equation}
d\equiv \frac{\sum_{i=k+1}^{p}\left( \mathbb{\tilde{V}}_{i,i}^{(1)}+\tilde{%
\mu}_{i}^{(1)}\cdot \tilde{\mu}_{i}^{(1)}\right) }{n}  \label{const}
\end{equation}%
to $T.$

\subsection{Corrections due to $3^{\text{rd}}$ and $4^{\text{th}}$-order
cumulants}

To find the contribution of the $D_{1}^{4}$ and $D_{1}^{6}$ terms, we repeat
the steps leading to (\ref{H2}), increasing the power of the $\frac{\partial 
}{\partial z_{1}}$ derivative (to $4$ and $6,$ respectively). This means that%
\begin{equation*}
\frac{\partial \exp (-\frac{z_{1}^{2}}{2})}{\partial z_{1}}=-z_{1}\exp (-%
\frac{z_{1}^{2}}{2})
\end{equation*}%
of the third last line of (\ref{H2}) needs to be replaced by%
\begin{equation*}
\frac{\partial ^{3}\exp (-\frac{z_{1}^{2}}{2})}{\partial z_{1}^{3}}%
=-(z_{1}^{3}-3z_{1})\exp (-\frac{z_{1}^{2}}{2})
\end{equation*}%
and by%
\begin{equation*}
\frac{\partial ^{5}\exp (-\frac{z_{1}^{2}}{2})}{\partial z_{1}^{5}}%
=-(z_{1}^{5}-10z_{1}^{3}+15z_{1})\exp (-\frac{z_{1}^{2}}{2})
\end{equation*}%
respectively, correspondingly modifying the rest.

For the $D_{1}^{4}$ contribution we now get, in place of (\ref{H2}):%
\begin{equation}
3\left( \frac{2u}{k}-\frac{2u^{2}}{k(k+2)}\right) \cdot \frac{u^{k/2-1}\exp
(-\frac{u}{2})}{2^{k/2}\Gamma (\frac{k}{2})}  \label{H4}
\end{equation}%
According to (\ref{pdf}), this needs to be multiplied by%
\begin{equation}
\frac{\sum_{i=1}^{k}\left( \tilde{\kappa}_{i,i,i,i}^{(4,0)}+4\tilde{\mu}%
_{i}^{(1)}\tilde{\kappa}_{i,i,i}^{(3,0)}\right) }{24n}  \label{k4}
\end{equation}%
before being added to (\ref{F}).

Similarly, $D_{1}^{6}$ leads to%
\begin{equation*}
15\left( -\frac{2u}{k}+\frac{4u^{2}}{k(k+2)}-\frac{2u^{3}}{k(k+2)(k+4)}%
\right) \cdot \frac{u^{k/2-1}\exp (-\frac{u}{2})}{2^{k/2}\Gamma (\frac{k}{2})%
}
\end{equation*}%
and is to be multiplied by%
\begin{equation}
\frac{\sum_{i=1}^{k}\left( \kappa _{i,i,i}^{(3,0)}\right) ^{2}}{72n}
\label{k6a}
\end{equation}

To deal with $D_{1}^{2}D_{2}^{2}$ (which covers the case of any $%
D_{i}^{2}D_{j}^{2}$ with $i\neq j$) we follow a similar procedure (skipping
the first rather obvious step):%
\begin{eqnarray*}
&&(2\pi )^{-k/2}\didotsint_{\mathcal{R}_{k}(\sqrt{u})}\frac{\partial ^{2}}{%
\partial z_{1}^{2}}\frac{\partial ^{2}}{\partial z_{2}^{2}}\exp \left( -%
\frac{\sum_{i=1}^{k}z_{i}^{2}}{2}\right) dz_{1}dz_{2}...dz_{k} \\
&=&(2\pi )^{-k/2}\dint\limits_{-\sqrt{u}}^{\sqrt{u}}\frac{\partial ^{2}\exp
(-\frac{z_{1}^{2}}{2})}{\partial z_{1}^{2}}\dint\limits_{-\sqrt{u-z_{1}^{2}}%
}^{\sqrt{u-z_{1}^{2}}}\frac{\partial ^{2}\exp (-\frac{z_{2}^{2}}{2})}{%
\partial z_{2}^{2}}\cdot \\
&&\didotsint_{\mathcal{R}_{k-2}(\sqrt{u-z_{1}^{2}-z_{2}^{2}})}\exp \left( -%
\frac{\sum_{i=3}^{k}z_{i}^{2}}{2}\right) dz_{3}...dz_{k}~dz_{2}~dz_{1} \\
&=&(2\pi )^{-k/2}\dint\limits_{-\sqrt{u}}^{\sqrt{u}}\frac{\partial ^{2}\exp
(-\frac{z_{1}^{2}}{2})}{\partial z_{1}^{2}}\dint\limits_{-\sqrt{u-z_{1}^{2}}%
}^{\sqrt{u-z_{1}^{2}}}\frac{\partial ^{2}\exp (-\frac{z_{2}^{2}}{2})}{%
\partial z_{2}^{2}}\cdot \\
&&\int_{0}^{\sqrt{u-z_{1}^{2}-z_{2}^{2}}}S_{k-2}(r)\exp \left( -\frac{r^{2}}{%
2}\right) dr~dz_{2}~dz_{1} \\
&=&-\frac{2\pi ^{-1}}{2^{k/2}\Gamma (\frac{k-2}{2})}\dint\limits_{-\sqrt{u}%
}^{\sqrt{u}}\frac{\partial ^{2}\exp (-\frac{z_{1}^{2}}{2})}{\partial
z_{1}^{2}}\dint\limits_{-\sqrt{u-z_{1}^{2}}}^{\sqrt{u-z_{1}^{2}}}\frac{%
\partial \exp (-\frac{z_{2}^{2}}{2})}{\partial z_{2}}\cdot \\
&&\frac{\partial }{\partial z_{2}}\int_{0}^{\sqrt{u-z_{1}^{2}-z_{2}^{2}}%
}r^{k-3}\exp \left( -\frac{r^{2}}{2}\right) dr~dz_{2}~dz_{1} \\
&=&-\frac{2\pi ^{-1}}{2^{k/2}\Gamma (\frac{k-2}{2})}\cdot \exp (-\frac{u}{2}%
)\dint\limits_{-\sqrt{u}}^{\sqrt{u}}(z_{1}^{2}-1)\dint\limits_{-\sqrt{%
u-z_{1}^{2}}}^{\sqrt{u-z_{1}^{2}}%
}z_{2}^{2}(u-z_{1}^{2}-z_{2}^{2})^{(k-4)/2}dz_{2}~dz_{1} \\
&=&-\frac{\pi ^{-1/2}}{2^{k/2}\Gamma (\frac{k+1}{2})}\cdot \exp (-\frac{u}{2}%
)\dint\limits_{-\sqrt{u}}^{\sqrt{u}%
}(z_{1}^{2}-1)(1-z_{1}^{2})^{(k-1)/2}dz_{1} \\
&=&\frac{2(1-\frac{u}{k+2})}{k\cdot 2^{k/2}\Gamma (\frac{k}{2})}\cdot
u^{k/2}\exp (-\frac{u}{2})=\left( \frac{2u}{k}-\frac{2u^{2}}{k(k+2)}\right)
\cdot \frac{u^{k/2-1}\exp (-\frac{u}{2})}{2^{k/2}\Gamma (\frac{k}{2})}
\end{eqnarray*}%
This is yet to be multiplied - see (\ref{pdf}) - by%
\begin{equation*}
\frac{\sum_{i\neq j}^{k}\left( \tilde{\kappa}_{i,i,j,j}^{(4,0)}+\tilde{\kappa%
}_{i,j,i,j}^{(4,0)}+\tilde{\kappa}_{i,j,j,i}^{(4,0)}+4\tilde{\mu}_{i}^{(1)}%
\tilde{\kappa}_{i,j,j}^{(3,0)}+4\tilde{\mu}_{i}^{(1)}\tilde{\kappa}%
_{j,i,j}^{(3,0)}+4\tilde{\mu}_{i}^{(1)}\tilde{\kappa}_{j,j,i}^{(3,0)}\right) 
}{24n}
\end{equation*}%
which, \emph{together} with (\ref{k4}) \emph{multiplied} by $3$ (how very
convenient of (\ref{H4}) to provide this extra factor) can be combined into%
\begin{equation}
\frac{\sum_{i,j=1}^{k}\left( \tilde{\kappa}_{i,i,j,j}^{(4,0)}+\tilde{\kappa}%
_{i,j,i,j}^{(4,0)}+\tilde{\kappa}_{i,j,j,i}^{(4,0)}+4\tilde{\mu}_{i}^{(1)}%
\tilde{\kappa}_{i,j,j}^{(3,0)}+4\tilde{\mu}_{i}^{(1)}\tilde{\kappa}%
_{j,i,j}^{(3,0)}+4\tilde{\mu}_{i}^{(1)}\tilde{\kappa}_{j,j,i}^{(3,0)}\right) 
}{24n}  \label{jv2}
\end{equation}%
Due to the total symmetry of cumulants, (\ref{jv2}) can be simplified
further to%
\begin{equation}
\frac{\sum_{i,j=1}^{k}\left( \tilde{\kappa}_{i,i,j,j}^{(4,0)}+4\tilde{\mu}%
_{i}^{(1)}\tilde{\kappa}_{i,j,j}^{(3,0)}\right) }{8n}  \label{jv3}
\end{equation}%
In (\ref{jv3}), we can now change the upper limit of the summation from $k$
to $p,$ as all those extra terms are equal to zero (as explained earlier).
And, since the resulting expression is rotationally invariant, we can
express it in terms of the old $X$ cumulants thus%
\begin{equation}
b\equiv \frac{\sum_{i,j=1}^{p}\left( \kappa _{i,i,j,j}^{(4,0)}+4\mu
_{i}^{(1)}\kappa _{i,j,j}^{(3,0)}\right) }{8n}  \label{b}
\end{equation}

Dealing with $D_{1}^{4}D_{2}^{2}$ is now quite easy - one has to replace (in
the previous derivation) 
\begin{equation*}
\frac{\partial ^{2}\exp (-\frac{z_{1}^{2}}{2})}{\partial z_{1}^{2}}%
=(z_{1}^{2}-1)\exp (-\frac{z_{1}^{2}}{2})
\end{equation*}%
by%
\begin{equation*}
\frac{\partial ^{4}\exp (-\frac{z_{1}^{2}}{2})}{\partial z_{1}^{4}}%
=(z_{1}^{4}-6z_{1}^{2}+3)\exp (-\frac{z_{1}^{2}}{2})
\end{equation*}%
and work out the details, getting%
\begin{equation*}
3\cdot \left( -\frac{2u}{k}+\frac{4u^{2}}{k(k+2)}-\frac{2u^{3}}{k(k+2)(k+4)}%
\right) \cdot \frac{u^{k/2-1}\exp (-\frac{u}{2})}{2^{k/2}\Gamma (\frac{k}{2})%
}
\end{equation*}%
yet to be multiplied by%
\begin{equation}
\frac{\sum_{i\neq j}^{k}\left( 
\begin{array}{c}
2\tilde{\kappa}_{i,i,j}^{(3,0)}\tilde{\kappa}_{j,j,j}^{(3,0)}+2\tilde{\kappa}%
_{i,j,i}^{(3,0)}\tilde{\kappa}_{j,j,j}^{(3,0)}+2\tilde{\kappa}%
_{j,i,i}^{(3,0)}\tilde{\kappa}_{j,j,j}^{(3,0)}+ \\ 
\left( \tilde{\kappa}_{i,j,j}^{(3,0)}\right) ^{2}+\left( \tilde{\kappa}%
_{j,i,j}^{(3,0)}\right) ^{2}+\left( \tilde{\kappa}_{j,j,i}^{(3,0)}\right)
^{2}+2\tilde{\kappa}_{i,j,j}^{(3,0)}\tilde{\kappa}_{j,i,j}^{(3,0)}+2\tilde{%
\kappa}_{i,j,j}^{(3,0)}\tilde{\kappa}_{j,j,i}^{(3,0)}+2\tilde{\kappa}%
_{j,i,j}^{(3,0)}\tilde{\kappa}_{j,j,i}^{(3,0)}%
\end{array}%
\right) }{72n}  \label{k6b}
\end{equation}

We would hope that by now the reader can supply the details of the $%
D_{1}^{2}D_{2}^{2}D_{3}^{2}$ derivation, getting the (at this point
expected) result of%
\begin{equation*}
\left( -\frac{2u}{k}+\frac{4u^{2}}{k(k+2)}-\frac{2u^{3}}{k(k+2)(k+4)}\right)
\cdot \frac{u^{k/2-1}\exp (-\frac{u}{2})}{2^{k/2}\Gamma (\frac{k}{2})}
\end{equation*}

Referring back to (\ref{pdf}), the last expression has to be multiplied by%
\begin{equation}
\frac{\sum_{i\neq j\neq \ell }^{k}\left( 
\begin{array}{c}
\tilde{\kappa}_{i,j,\ell }^{(3,0)}\cdot (\tilde{\kappa}_{i,j,\ell }^{(3,0)}+%
\tilde{\kappa}_{i,\ell ,j}^{(3,0)}+\tilde{\kappa}_{j,i,\ell }^{(3,0)}+\tilde{%
\kappa}_{j,\ell ,i}^{(3,0)}+\tilde{\kappa}_{\ell ,i,j}^{(3,0)}+\tilde{\kappa}%
_{\ell ,j,i}^{(3,0)})+ \\ 
\tilde{\kappa}_{i,j,j}^{(3,0)}\cdot \tilde{\kappa}_{i,\ell ,\ell }^{(3,0)}+%
\tilde{\kappa}_{i,j,j}^{(3,0)}\cdot \tilde{\kappa}_{\ell ,i,\ell }^{(3,0)}+%
\tilde{\kappa}_{i,j,j}^{(3,0)}\cdot \tilde{\kappa}_{\ell ,\ell ,i}^{(3,0)}+%
\tilde{\kappa}_{j,i,j}^{(3,0)}\cdot \tilde{\kappa}_{i,\ell ,\ell }^{(3,0)}+
\\ 
\tilde{\kappa}_{j,i,j}^{(3,0)}\cdot \tilde{\kappa}_{\ell ,i,\ell }^{(3,0)}+%
\tilde{\kappa}_{j,i,j}^{(3,0)}\cdot \tilde{\kappa}_{\ell ,\ell ,i}^{(3,0)}+%
\tilde{\kappa}_{j,j,i}^{(3,0)}\cdot \tilde{\kappa}_{i,\ell ,\ell }^{(3,0)}+%
\tilde{\kappa}_{j,j,i}^{(3,0)}\cdot \tilde{\kappa}_{\ell ,i,\ell }^{(3,0)}+
\\ 
\tilde{\kappa}_{j,j,i}^{(3,0)}\cdot \tilde{\kappa}_{\ell ,\ell ,i}^{(3,0)}%
\end{array}%
\right) }{72n}  \label{k6}
\end{equation}%
which again, rather conveniently, can be combined with (\ref{k6a}) and (\ref%
{k6b}) since%
\begin{equation*}
\sum_{i,j,\ell =1}^{k}...=\sum_{i\neq j\neq \ell }^{k}...+\sum_{i=j\neq \ell
}^{k}...+\sum_{i\neq j=\ell }^{k}...+\sum_{i=\ell \neq
j}^{k}...+\sum_{i=j=\ell }^{k}...
\end{equation*}%
resulting in%
\begin{equation}
\sum_{i,j,\ell =1}^{k}\left( -"-\right) =\sum_{i\neq j\neq \ell }^{k}\left(
-"-\right) +3\sum_{i=j}^{k}\left( -""-\right) +15\sum_{i}^{k}\left( \kappa
_{i,i,i}^{(3,0)}\right) ^{2}  \label{jv4}
\end{equation}%
where $(-"-)$ refers to the big parentheses of (\ref{k6}) and $(-""-)$ to
the big parentheses of (\ref{k6b}). Due to the total symmetry of the
cumulants, the full answer can be further simplified to 
\begin{equation*}
\sum_{i,j,\ell =1}^{k}\frac{\left( \tilde{\kappa}_{i,j,\ell }^{(3,0)}\right)
^{2}}{12n}+\sum_{i,j,\ell =1}^{k}\frac{\tilde{\kappa}_{i,j,j}^{(3,0)}\cdot 
\tilde{\kappa}_{i,\ell ,\ell }^{(3,0)}}{8n}
\end{equation*}%
As before, increasing the upper limit of the summation from $k$ to $p$
changes nothing. And, since both terms are rotationally symmetric, the final
version of the formula is 
\begin{equation}
c\equiv \sum_{i,j,\ell =1}^{p}\frac{\left( \kappa _{i,j,\ell
}^{(3,0)}\right) ^{2}}{12n}+\sum_{i,j,\ell =1}^{p}\frac{\kappa
_{i,j,j}^{(3,0)}\cdot \kappa _{i,\ell ,\ell }^{(3,0)}}{8n}  \label{c}
\end{equation}%
using the original $\mathbf{X}$ (not the rotated $\mathbf{Z}$) cumulants.

\section{Conclusion}

We have thus found that the $\frac{1}{n}$-accurate CDF of $T-d$ is given by%
\begin{eqnarray*}
&&F(u)+\left[ a\cdot \left( -\frac{2u}{k}\right) +b\cdot \left( \frac{2u}{k}-%
\frac{2u^{2}}{k(k+2)}\right) +\right. \\
&&\left. c\cdot \left( -\frac{2u}{k}+\frac{4u^{2}}{k(k+2)}-\frac{2u^{3}}{%
k(k+2)(k+4)}\right) \right] \cdot \frac{u^{k/2-1}\exp (-\frac{u}{2})}{%
2^{k/2}\Gamma (\frac{k}{2})}
\end{eqnarray*}%
where $a,$ $b,$ $c$ and $d$ are defined in (\ref{a}), (\ref{b}), (\ref{c})
and (\ref{const}) respectively, and $F(u)$ is the CDF of the $\chi _{k}^{2}$
distribution defined in (\ref{F}).

The corresponding PDF equals, by simple differentiation%
\begin{eqnarray*}
&&\frac{u^{k/2-1}\exp (-\frac{u}{2})}{2^{k/2}\Gamma (\frac{k}{2})}\cdot %
\left[ 1+a\cdot \left( \frac{u}{k}-1\right) +b\cdot \left( \frac{u^{2}}{%
k(k+2)}-\frac{2u}{k}+1\right) +\right. \\
&&\left. c\cdot \left( \frac{u^{3}}{k(k+2)(k+4)}-\frac{3u^{2}}{k(k+2)}+\frac{%
3u}{k}-1\right) \right]
\end{eqnarray*}%
when $u>0$ (zero otherwise). Note that the actual diagonalization of $%
\mathbb{V}_{0}$ is needed only for establishing the $a$ and $d$ constants - $%
b$ and $c$ can be found directly from the $3^{rd}$ and $4^{th}$ cumulants of 
$\mathbf{X}.$

Examples of applying this procedure can be found in our two references; the
second one incorrectly assumed that $d=0$ - with the help of this article,
one can easily make the corresponding correction (it turns out that, in this
particular case, $a$ and $d$ must have equal values).

And one final remark: it would prove rather difficult to extend this
procedure to the full $\frac{1}{n^{2}}$ accuracy; nevertheless, it usually
proves beneficial to compute the value of each $a$ and $d$ to be $\frac{1}{%
n^{2}}$ accurate - this tend to extend the applicability of the resulting
approximation to such small sample sizes that further improvements do not
appear (from a practical point of view) necessary.

\end{document}